\documentclass[a4paper,10pt]{article}
 \usepackage{hyperref}
\usepackage{paper-en}

\def\thetitle{Milnor's cartography problem}
\def\theauthors{Anton Petrunin}

\hypersetup{colorlinks=true,
citecolor=black,
linkcolor=black,
anchorcolor=black,
filecolor=black,
menucolor=black,
urlcolor=black,
pdftitle={\thetitle},
pdfauthor={\theauthors}
}

\begin{document}
\title{\thetitle}
\author{\theauthors}
\date{}
\maketitle

\begin{abstract}
We solve John Milnor's problem: among all convex regions of a given area on the sphere, the round disk requires the greatest distortion in cartographic projections onto the plane.
\end{abstract}

\section{Introduction}

Let $\mathbb{S}^2$ be the unit sphere, equipped with the geodesic metric, and let $\Omega\subset\mathbb{S}^2$.
Recall that a map $f\:\Omega\to\RR^2$ is called \emph{bi-Lipschitz} if
\[L_1\cdot|x-y|_{\mathbb{S}^2}\le |f(x)-f(y)|_{\RR^2}\le L_2\cdot|x-y|_{\mathbb{S}^2}\]
for some positive constants $L_1$ and $L_2$ and any two points $x,y\in \Omega$;
here $|\ -\ |_{\mathbb{S}^2}$ and $|\ -\ |_{\RR^2}$ stand for the geodesic distance on the sphere $\mathbb{S}^2$ and on the Euclidean plane $\RR^2$, respectively.

Given a set $\Omega$, we want to construct its cartographic projection with small distortion ratio $L_2/L_1$.
After rescaling, we may assume that $L_1=1$ (so $f$ is noncontracting) and minimize $L_2$.

Let us denote by $D_\rho$ the round disk of radius $\rho$ in $\mathbb{S}^2$.
It is straightforward to check that the polar coordinates with the origin at the center of $D_\rho$ define a
map $f\:D_\rho\to\RR^2$ with the distortion ratio $\tfrac{\rho}{\sin \rho}$;
in cartography, this is known as the \emph{azimuthal equidistant projection}.

John Milnor \cite{milnor} observed that the value $\tfrac{\rho}{\sin \rho}$ is optimal;
that is, \textit{if $f\:D_\rho\z\to\RR^2$ is a bi-Lipschitz map,
then its distortion ratio must be at least $\tfrac{\rho}{\sin \rho}$; moreover, if equality holds, then, up to a similarity of the plane, the map $f$ is the same as above.}
In the same paper, he formulated the following problem; see also the expository article by Étienne Ghys \cite{ghys}.

\begin{thm}{Problem}\label{prob:main}
Among all convex regions in $\mathbb{S}^2$ of a given area, the disk requires the largest distortion ratio.

More precisely, any closed convex set $\Omega\subset\mathbb{S}^2$ such that $\area \Omega=\area D_\rho$ admits a map $\Omega\to\RR^2$ with Lipschitz constants $1$ and $\tfrac \rho{\sin \rho}$.
\end{thm}

The proof proceeds by subdividing $\Omega$ into thin slices such that the union of their borderlines forms a geodesic tree.
This tree is developed in the plane, preserving edge lengths and angles.
Then the resulting map on the tree is extended to each slice.
We impose a certain balancing condition on the subdivision, which makes it possible to keep the distortion under control.

Our solution was found while writing a note \cite{petrunin2026} for school students about Gary Lawlor's proof of the isoperimetric inequality \cite{lawlor}.
In fact, his construction produces a subdivision that can be used to solve a version of Milnor's problem with perimeter instead of area; this has been known for a while \cite{petrunin2016,petrunin2025}.

\section{Balanced sets}

\begin{wrapfigure}{o}{36mm}
\centering
\vskip-12mm
\includegraphics{mppics/pic-10}
\vskip-0mm
\end{wrapfigure}

Fix a closed convex set $\Omega\subset\mathbb{S}^2$.
A convex subset $K\z\subset \Omega$ is called \emph{admissible} if it is bounded by a (possibly degenerate) arc of $\partial \Omega$ and a polygonal line $a_0a_1\dots a_n$ joining the endpoints of this arc.
The arc $\delta=\partial \Omega\cap K$ will be called the \emph{coastline} of $K$;
the polygonal line $a_0a_1\dots a_n$ will be called the \emph{borderline} of $K$;
their common endpoints $a_0$ and $a_n$ will be called the \emph{corners} of $K$.
Let us define the \emph{turn} of $K$ (briefly, $\turn K$) as the sum of exterior angles of $K$ at $a_1$, $a_2,\dots,a_{n-1}$.

Fix $\rho\le \tfrac\pi2$ and define
\[\balance_\rho K\df\lambda(\rho)\cdot\area K+\turn K-\pi,\]
where
\[\lambda(\rho)\df\frac{\area D_{\frac\pi2}}{\area D_\rho}=\frac{1}{1-\cos \rho}.\]
Since $\rho\le \tfrac\pi2$, we have $\lambda(\rho)\ge 1$.
If $\balance_\rho K=0$, then $K$ is called \emph{$\rho$-balanced}.
Note that every sector of $D_\rho$ is $\rho$-balanced.

\begin{thm}{Lemma}\label{lem:dist-r}
Any point of a $\rho$-balanced subset $K$ lies at distance at most $\rho$ from its coastline.
Moreover, if the maximum distance is $\rho$, then $K$ is congruent to a sector of $D_\rho$.
\end{thm}

Recall that for any spherical $n$-gon $P$ we have
\[\theta_1+\ldots+\theta_n=2\cdot\pi-\area P,\]
where $\theta_1,\ldots,\theta_n$ denote the signed external angles of $P$; thus $\theta_i\ge 0$ if the vertex is convex, $\theta_i\le 0$ if the vertex is concave, and $|\theta_i|<\pi$ for any $i$.
This is the \emph{spherical area formula}, which is a special case of the Gauss--Bonnet theorem.

\parit{Proof.}
Assume there is a point $x\in K$ at distance at least $\rho$ from the coastline.
Let $L$ be the curvilinear triangle with vertex $x$ and the coastline as the base.
Denote by $S$ the intersection of $L$ with the $\rho$-disk centered at $x$.

\begin{wrapfigure}{o}{36mm}
\centering
\vskip-1mm
\includegraphics{mppics/pic-20}
\vskip-0mm
\end{wrapfigure}

Note that $S$ is a sector of the $\rho$-disk.
Moreover, $S$ is convex; otherwise $\area K>\tfrac12\cdot\area D_\rho$, so $K$ is not $\rho$-balanced.

It follows that $S$ is $\rho$-balanced.
Furthermore, $\turn S=\turn L$, and since $S\subset L$, we have $\area S\z\le \area L$.
Therefore,
\[0=\balance_\rho S\le \balance_\rho L;\]
moreover, if equality holds, then $L=S$.

On the other hand, $L\subset K$, and so $\area L\z\le \area K$.
Furthermore, applying the spherical area formula to the complement $K\setminus L$, we get
\[\turn L-\turn K\le\area K-\area L.\]
Since $\lambda(\rho)\ge1$, we have
\[0=\balance_\rho K\ge \balance_\rho L;\]
moreover, in the case of equality, we have $K=L$, which finishes the proof.
\qeds

Applying the spherical area formula, we also obtain the following.

\begin{wrapfigure}{o}{28mm}
\centering
\vskip-4mm
\includegraphics{mppics/pic-25}
\vskip-2mm
\end{wrapfigure}

\begin{thm}{Observation}\label{lem:nondegenerate}
Let $K$ be an admissible set.
Suppose $\balance_\rho K\le 0$.
Then the coastline of $K$ is nondegenerate; in other words, $K$ has two distinct corners, say $l$ and~$r$.
Furthermore, if $\phi$ and $\psi$ denote the angles that the chord $[l,r]$ makes with the borderline, then $\phi+\psi\le \pi$.
\end{thm}

Let $K$ be an admissible set and $S$ be a sector of a spherical disk.
We say that $K$ is \emph{properly inscribed} in $S$ if the corners of $K$ lie on two different sides of $S$ and the coastline of $K$ cuts off from $S$ a convex set that contains $K$ and the center of $S$.

{

\begin{wrapfigure}{o}{36mm}
\centering
\vskip-4mm
\includegraphics{mppics/pic-30}
\vskip-2mm
\end{wrapfigure}

\begin{thm}{Lemma}\label{lem:angle}
If an admissible set $K$ is properly inscribed in a sector $S$ of $D_\rho$,
then $\balance_\rho K\le 0$ and, in the case of equality, we have $K=S$.
\end{thm}

\parit{Proof.}
Let us denote by $L$ the convex subset that is cut off from $S$ by the coastline of $K$.
Arguing as in \ref{lem:dist-r}, we get
\[0=\balance_\rho S\ge \balance_\rho L\ge \balance_\rho K.\]
If $\balance_\rho K=0$, then we have two equalities, and hence $S=L\z=K$.
\qeds

}

\section{Slicing}

\begin{thm}{Proposition}\label{prop:subdivision}
Let $\Omega\subset \mathbb{S}^2$ be a closed convex set.
If $\area \Omega=\area D_\rho$, then $\Omega$ admits a finite
subdivision into $\rho$-balanced sets with arbitrarily short coastlines.
\end{thm}

\begin{thm}{Observation}\label{obs:subdivision}
Suppose that an admissible set $K$ is subdivided into two admissible sets $L$ and $R$ by a chord from the borderline to the coastline.
Then
\[
\balance_\rho K=\balance_\rho L+\balance_\rho R.
\]

\end{thm}

\parit{Proof of \ref{prop:subdivision}.}
Let us split $\Omega$ by a chord into two sets $K_1$ and $K_2$ of equal area.
Since the chord has zero turn, each $K_i$ is $\rho$-balanced.

\begin{wrapfigure}{o}{30mm}
\centering
\vskip-4mm
\includegraphics{mppics/pic-50}

\end{wrapfigure}

It remains to show that any $\rho$-balanced subset $K$ can be subdivided into a finite number of $\rho$-balanced sets with arbitrarily short coastlines.

Let $l$ and $r$ be the left and right corners of $K$.
Choose a point $p$ on the borderline and a point $q$ on the coastline.
Then the chord $[p,q]$ splits $K$ into two admissible sets $L\ni l$ and $R\ni r$.

If $p$ and $q$ are close to $l$, then by \ref{lem:angle}, $\balance_\rho L\z\le 0$, and, by the observation, $\balance_\rho R\z\ge0$.
Similarly, if $p$ and $q$ are close to $r$, then $\balance_\rho R\le 0$, and $\balance_\rho L\ge0$.
Since $\balance_\rho L$ depends continuously on the choice of the chord, we can choose an intermediate chord $[p,q]$ that splits $K$ into two $\rho$-balanced sets.

Furthermore, for any $\eps>0$ we can choose $\delta>0$ such that if the coastline of $K$ is longer than $\eps$, then
the position of $q$ in the above construction can be chosen $\delta$-far from the corners along the coastline.
Indeed, we can assume that $\eps<\rho$.
Choose $\delta=\tfrac\eps2$.
Let $p$ be a point on the first edge of the borderline at distance less than $\tfrac\rho2$ from $l$.
If $q$ is $\delta$-far from $l$ along the coastline, then \ref{lem:angle} implies $\balance_\rho L\z\le 0$.

It follows that the coastline of each of the slices $L$ and $R$ is shorter than that of $K$ by at least $\delta$.
Therefore, repeating this construction recursively yields the needed subdivision.
\qeds

\section{Chord foliation}

\begin{thm}{Proposition}\label{prop:diam-r}
Let $K$ be an admissible set.
Suppose that $\balance_\rho K\le 0$ and every chord from the borderline to the coastline has length less than $\rho$.
Then the interior of $K$ can be foliated by chords from the borderline to the coastline such that the slice between any pair of chords, or on one side of a chord, has nonpositive $\rho$-balance.
\end{thm}

\parit{Proof.}
Let $[p,q]$ be a chord of $K$ from the borderline to the coastline.
Denote by $L_{[p,q]}$ the part of $K$ on the left side of $[p,q]$.

Choose an orientation of both the borderline and the coastline of $K$ from the left corner to the right one.
Note that the function
\[(p,q)\mapsto \balance_\rho L_{[p,q]}\]
is increasing in $p$ and decreasing in $q$ with respect to the orientation.

Indeed, by moving $p$ slightly along the orientation, we increase the area and turn of $L_{[p,q]}$;
the latter follows since in a spherical triangle of diameter smaller than $\tfrac\pi2$, an exterior angle is greater than an interior angle at another vertex (recall that $\rho\le \tfrac\pi2$).
Further, when $q$ moves along the orientation, the area increases and the turn decreases, but the condition on chord length ensures that the latter term wins.

\begin{wrapfigure}{o}{30mm}
\centering
\vskip-6mm
\includegraphics{mppics/pic-60}
\vskip0mm
\end{wrapfigure}

It follows that the chords $[p,q]$ defined by the equation
\[\balance_\rho L_{[p,q]}=0\]
do not intersect and foliate a subset $K'\subset K$.

If $K'\ne K$ then the foliation has a chord with an endpoint at the right corner.
By \ref{obs:subdivision}, this chord bounds an admissible subset with nonpositive $\rho$-balance.
By \ref{lem:nondegenerate}, its coastline cannot be degenerate.
Therefore, the other endpoint of the chord must be on the coastline.
In this case, the foliation can be extended by chords with an endpoint at this corner.
By the condition on chord length, the extended foliation meets the conditions.
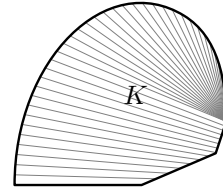
\qeds

\section{Fishbone map}

Let $\beta$ be the borderline of an admissible set $K$, and let $\tilde \beta$ be its \emph{development} in the plane;
that is, $\tilde \beta$ is a polygonal line in $\RR^2$ with the same side lengths and oriented angles as $\beta$.
Note that there is a natural length-preserving map $\beta\to\tilde \beta$.

Suppose $\balance_\rho K\le 0$; in particular $\turn K<\pi$.
Then $\tilde \beta$ and extensions of its first and last segments bound a convex set of $\RR^2$; let us call it the \emph{plane slice} of $K$ and denote it by $\tilde K$.

Suppose the interior of $K$ is foliated by chords from the borderline to the coastline.
For each chord $\ell=[p,q]$ in the foliation, consider the corresponding half-line $\tilde \ell$ in $\tilde K$;
it starts at the point $\tilde p\in\tilde \beta$ that corresponds to $p\in\beta$ and makes the same angle with $\tilde \beta$.
Suppose that the original foliation is \emph{spreading}; that is, every slice between two chords of the foliation, or on one side of a chord, has turn at most~$\pi$.
Then the resulting half-lines are pairwise disjoint and they foliate a region in the interior of $\tilde K$.
\begin{figure}[ht!]
\centering
\vskip-0mm
\includegraphics{mppics/pic-70}
\end{figure}
There might be one or two angles in $\tilde K$ adjacent to its infinite sides that are not covered by this foliation.
In this case, foliate them by the half-lines starting from the vertices of these angles; these vertices correspond to the corners of $K$.

The natural map $\beta\to\tilde\beta$ extends to the \emph{fishbone map} $f\:K\to \tilde K$ that sends each chord of the foliation isometrically to the starting segment of the corresponding half-line.

\begin{thm}{Observation}\label{obs:fishbone}
Let $K$ be an admissible set such that $\balance_\rho K\le 0$ and every chord from the borderline to the coastline has length less than $\rho$.
Then the foliation of $K$ provided by \ref{prop:diam-r} induces a $\tfrac \rho{\sin \rho}$-Lipschitz fishbone map $f\:K\to \tilde K$.
\end{thm}

\parit{Proof.}
Choose an oriented orthonormal frame $\vec v$, $\vec w$ in the interior of $K$ such that $\vec v$ points in the direction of the chord foliation.
Further, choose an oriented orthonormal frame $\tilde {\vec v}$, $\tilde {\vec w}$ in the interior of the plane slice $\tilde K$ such that $\tilde{\vec v}$ points in the direction of the half-line foliation.

Let $\ell_n\to\ell$ be a converging sequence of distinct chords in the foliation; let $s$ be the arc-length parameter of $\ell$ starting at its borderline point.
After rescaling and passing to a subsequence of $\ell_n$, the normal components of the variations converge to a Jacobi field, say $\vec i$; we can assume that
\begin{align*}
\vec i(s)&=\vec w(s)\cdot\sin(a+s)
\intertext{for some $a$.
The corresponding Jacobi field in the plane is}
\tilde{\vec i}(s)&=\tilde{\vec w}(s)\cdot (\sin a+s\cdot \cos a).
\intertext{Let}
\nu(s)&\df\frac{|\tilde{\vec i}(s)|}{|\vec i(s)|}=\frac{\sin a+s\cdot \cos a}{\sin(a+s)}.
\end{align*}

It is sufficient to show that
\[\nu(s)\le \frac \rho{\sin \rho}.\leqno({*})\]
Indeed, since the sequence $\ell_n$ was arbitrary, the expansion of $f$ in the direction of $\vec w$ is at most $\tfrac \rho{\sin \rho}$.
The expansion in the direction of $\vec v$ is $1$.
By the Gauss lemma, $f$ is locally $\tfrac \rho{\sin \rho}$-Lipschitz.
Since $K$ is convex, it is globally $\tfrac \rho{\sin \rho}$-Lipschitz.

Applying the balance inequality to the slices between the converging chords, we get its infinitesimal version:
\[\lambda(\rho)\cdot[\cos a- \cos (a+s)]=\lambda(\rho)\cdot\int\limits_0^{s}\sin (a+t)\cdot dt\le\cos a;
\leqno({*}{*})\]
here the integral of $\sin (a+t)$ plays the role of $\area K$ and $\cos a=\sin' a$ plays the role of $\pi-\turn K$.

Note that $0\le a\le a+s\le \tfrac\pi2$; the last inequality follows from \ref{lem:nondegenerate}.
Let \[\mu=\cos a-\cos (a+s);\] by $({*}{*})$, we have $\mu\le 1-\cos\rho$.

Clearly, $\nu$ is nondecreasing in $s$ for fixed $a$.
The computations below show that $\nu$ is nonincreasing in $a$ for fixed $\mu$.
Thus $\nu$ is maximized when $a=0$ and $s=\rho$, where $\nu=\tfrac{\rho}{\sin\rho}$; hence $({*})$ follows.

\parit{Computations.}
Regard $s$ and $\nu$ as functions of $a$.
For fixed $\mu$, we have $(a+s)'\z=\tfrac{\sin a}{\sin(a+s)}$.
Since $0\le a\le a+s\le \tfrac\pi2$, we have $\sin s\le s$ and $\mu\cdot\cos (a+s)\ge 0$.
Hence
\[\nu'=\frac{\sin a\cdot (\sin s -s\cdot [\mu\cdot\cos(a+s)+1])}{\sin^3 (a+s)}\le0.\]
\qedsf

The fishbone map $f$ has a natural left inverse $b\:\tilde K\to \mathbb{S}^2$; let us call it the \emph{backbone map}.
It maps the half-line $\tilde \ell$ in the foliation of $\tilde K$ to the geodesic that extends the corresponding chord $\ell$ in a locally distance-preserving way.
The additional half-lines in the foliation of $\tilde K$ are mapped to geodesics making the same angle with $\beta$ and starting at the corresponding corners of $K$.

\begin{thm}{Observation}\label{obs:backbone}
Let $K$ be an admissible set such that $\balance_\rho K\le 0$ and every chord from the borderline to the coastline has length less than $\rho$.
Then the foliation of $K$ provided by \ref{prop:diam-r} induces a $1$-Lipschitz backbone map.
\end{thm}

\parit{Proof.}
As before, we need to compare the absolute values of the Jacobi fields
\[\tilde {\vec i}(s)=\tilde {\vec w}(s)\cdot (\sin a+ s\cdot\cos a)\quad\text{and}\quad{\vec i}(s)= {\vec w}(s)\cdot \sin (a+s).\]
Recall that $0\le a\le \tfrac\pi2$;
therefore
\[|\tilde {\vec i}(s)|=|\sin a+ s\cdot \cos a|\ge |\sin (a+s)|=|{\vec i}(s)|\]
for any $s\ge0$, and the result follows.
\qeds

\section{Assembling the proof}

\parit{Proof of \ref{prob:main}.}
Since $\Omega$ is convex, we have $\rho\le \tfrac\pi2$.
Moreover, we can assume that $\rho<\tfrac\pi2$;
otherwise, $\Omega$ is a hemisphere and polar coordinates provide the required map.

Applying \ref{prop:subdivision}, we get a subdivision of $\Omega$ into $\rho$-balanced sets $K_1,\dots,K_m$.
The union of the borderlines of $K_i$ is a geodesic tree; denote it by $T$.

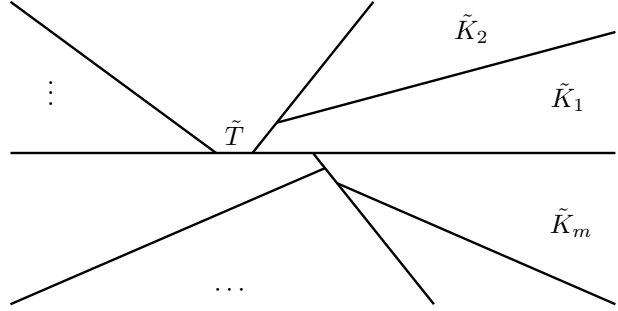
\begin{figure}[ht!]
\centering
\vskip-0mm
\includegraphics{mppics/pic-80}
\end{figure}

Consider the development of $T$;
it is a geodesic tree in the plane with the same edge lengths and the same oriented angles as in $T$.
Extend the terminal edges of the development beyond the terminal endpoints; denote the resulting tree by $\tilde T$.
Note that $\tilde T$ subdivides the plane into plane slices $\tilde K_1,\dots,\tilde K_m$ of $K_1,\dots,K_m$, respectively.

By \ref{prop:subdivision}, we can assume that the coastlines of $K_i$ are shorter than a given $\eps>0$.
Suppose $\eps$ is small; in particular, $\rho+\eps<\tfrac\pi2$.
We can apply \ref{prop:diam-r} with radius $\rho+\eps$ to each $K_i$.
Indeed, $\balance_{\rho+\eps}K_i<0$ and \ref{lem:dist-r} implies the condition on chord lengths.

The fishbone maps $f_i\: K_i\to \tilde K_i$ agree on $T$ and therefore fit together in one map $f\:\Omega\to\RR^2$.
By \ref{obs:fishbone}, $f$ is $\tfrac{\rho+\eps}{\sin(\rho+\eps)}$-Lipschitz.

Furthermore, the backbone maps $b_i\: \tilde K_i\to \mathbb{S}^2$ also agree on $\tilde T$, and they fit together in one map $b\:\RR^2\to\mathbb{S}^2$, which is $1$-Lipschitz by \ref{obs:backbone}.

By the construction, $b\circ f=\id_\Omega$.
Thus $f\:\Omega\to\RR^2$ has Lipschitz constants $1$ and $\tfrac{\rho+\eps}{\sin(\rho+\eps)}$.
Since $\eps>0$ can be chosen arbitrarily, the Arzelà--Ascoli theorem provides a map $f_0\:\Omega\to\RR^2$ with Lipschitz constants $1$ and $\tfrac\rho{\sin(\rho)}$.
\qeds

\section{Remarks}

The argument above suggests the following exercise.

\begin{thm}{Exercise}
Show that if $\Omega$ is not a round disk, then there is a map $\Omega\to\RR^2$ with distortion ratio strictly less than $\tfrac\rho{\sin \rho}$.
\end{thm}

Already in Milnor's paper, it was noted that the Arzelà--Ascoli theorem proves the existence of a map $f_0\:\Omega\to\RR^2$ with optimal Lipschitz constants.
We can assume that $f_0$ has lower Lipschitz constant $1$; therefore the upper Lipschitz constant $L_2$ is minimal.
Assuming that $\Omega$ is convex, we expect that $f_0$ has the following properties:
\begin{enumerate}[(a)]
\item There is a closed subset $X$ of $\Omega$ such that the differential $d_x f_0$ is an isometry at every $x\in X$, and the complement $\Omega\setminus X$ has a foliation by geodesics such that $f_0$ isometrically maps each geodesic in the foliation to a geodesic in the plane.
\item The local upper Lipschitz constants tend to $L_2$ as one approaches the boundary of $\Omega$.
\item The image $f_0(\Omega)$ is convex.
\end{enumerate}

\medskip

The following exercise addresses a related problem for nonconvex domains,
for example, the contiguous territory of your favorite country.
I know a proof, so I cannot honestly call this an open problem,
but if you solve it, then you should publish the solution.

\begin{thm}{Very advanced exercise}
Let $S$ be a noncomplete simply connected surface with Gauss curvature $1$ and area equal to that of $D_\rho$.
Show that there is a bi-Lipschitz homeomorphism from $S$ to a flat surface with distortion ratio at most~$\tfrac{\rho}{\sin\rho}$.
\end{thm}

The exercise provides a cartographic projection that can be made of paper, altho it may have overlaps when unfolded on a table.

In every map of the contiguous US that I know, the Lipschitz constants differ by more than $2\%$.
The exercise gives a map for which they differ by about $1\%$.
The actual optimal map should be even better, and there is a good chance that it can be unfolded without overlaps.

\paragraph{Acknowledgments.}
I want to thank
Arseniy Akopyan,
Roman Karasev,
Dmitri Panov,
and
Fedor Petrov
for their help.

{\sloppy
\def\emph{\textit}
\printbibliography[heading=bibintoc]\fussy}
\end{document}